\documentclass[a4paper, 12pt]{article}
\usepackage{amsmath}
\usepackage{amsfonts}
\usepackage{amssymb}
\usepackage{amscd}
\usepackage{amsthm}

\newcommand{\A}{{\mathbf A}}
\newcommand{\Df}{{\mathbf D}}
\newcommand{\Z}{{\mathbb Z}}

\newcommand{\Lc}{{\mathcal L}}
\newcommand{\Pc}{{\mathcal P}}
\newcommand{\M}{{\mathcal M}}

\newcommand{\Oc}{{\mathcal O}}

\newcommand{\m}{\mathfrak{m}}

\newcommand{\ord}{{\rm ord}}

\newcommand{\Pic}{{\rm Pic}}
\newcommand{\Hom}{{\rm Hom}}

\newcommand{\Nm}{{\rm Nm}}
\newcommand{\Div}{{\rm Div}}

\renewcommand{\div}{{\rm div}}

\newcommand{\Ad}{{\rm Ad}}

\theoremstyle{theorem}
\newtheorem{theor}{Theorem}[section]
\theoremstyle{theorem}
\newtheorem{prop}{Proposition}[section]
\theoremstyle{theorem}
\newtheorem{corol}{Corollary}[section]
\theoremstyle{theorem}
\newtheorem{lemma}{Lemma}[section]

\theoremstyle{definition}

\theoremstyle{remark}
\newtheorem{rmk}{Remark}[section]
\theoremstyle{remark}

\theoremstyle{remark}

\title{Poincar\'e biextension and ideles on an algebraic
curve}
\author{Sergey Gorchinskiy\footnote{The author was
partially supported by RFFI grants 04-01-00613 and
05-01-00455.}}
\date{}

\begin{document}
\maketitle
%\tableofcontents
%\bigskip
\vspace{4ex}

\begin{abstract}
Arbarello, de Concini, and Kac have constructed a central extension
of the ideles group on a smooth projective algebraic curve $C$. We
show that this central extension induces the theta-bundle on the
class group of degree $g-1$ divisors on $C$, where $g$ is the genus
of the curve $C$. The other result of the paper is the relation
between the product of the norms of the tame symbols over all points
of the curve, considered as a pairing on the ideles group, and the
Poincar\'e biextension of the Jacobian of $C$. As an application we
get a new proof of the adelic formula for the Weil pairing.
\end{abstract}

\section*{Introduction}

There exists a general ideology which tells that many notions and
constructions in algebraic geometry can be translated into the
language of certain adelic groups defined for a scheme and some
additional data on it, for instance a coherent sheaf (see more
details in \cite{FP} or \cite{Par}). This article provides a new
example to this approach.

Let $C$ be a smooth projective curve over a field $k$, $m$ be an
integer prime to ${\rm char}(k)$. Consider two divisors $D$ and $E$
on $C$ such that their classes in $\Pic(C)$ belong to the
$m$-torsion. Let $\alpha,\beta\in\A^*_C$ be two ideles such that
$\div(\alpha)=D$, $\div(\beta)=E$, and the ideles $\alpha^m$ and
$\beta^m$ are principal, i.e., belong to the subgroup
$k(C)^*\subset\A^*_C$; then the Weil pairing $\phi_m([D],[E])$ of
the classes of $E$ and $D$ in $\Pic(C)_m$ can be given by the
following adelic formula:
$$
\phi_m([D],[E])=(\prod_{x\in C}\Nm_{k(x)/k}[
(-1)^{\ord_x(\alpha_x)\ord_x(\beta_x)} (\alpha_x^{\ord_x(\beta_x)}
\beta_x^{-\ord_x(\alpha_x)})(x)])^m.
$$
The first proof of this formula when $C$ is of any genus appeared in
\cite{How}; a more elementary proof was given later by M. Mazo in
\cite{Maz}.

On the other hand, Arbarello, de Concini, and Kac have constructed
in \cite{AKC} a certain central extension of the ideles group
$$
0\to k^*\to\widehat{\A}^*_C \stackrel{\pi}\longrightarrow\A_C^*\to
0.
$$
It follows that the commutator in this extension is given by the
formula
$$
[\pi^{-1}(\alpha),\pi^{-1}(\beta)]=(-1)^{\deg(\alpha)\deg(\beta)}
\prod_{x\in C}\Nm_{k(x)/k}[ (-1)^{\ord_x(\alpha_x)\ord_x(\beta_x)}
(\alpha_x^{\ord_x(\beta_x)} \beta_x^{-\ord_x(\alpha_x)})(x)]
$$
for any ideles $\alpha,\beta\in\A_C^*$. The goal of this paper is to
find a reason for the apparent similarity of these two formulas. In
turns out that there is a close relation between the central
extension from \cite{AKC} and the Poincar\'e biextension over the
Jacobian of $C$, which defines the Weil pairing. Namely, we prove
that the Poincar\'e biextension is isomorphic to some quotient of
the canonically trivial biextension $\Lambda(\widehat{\A}^*_C)$
associated in a usual way with the central extension
$\widehat{\A}^*_C$.

The paper is organized as follows. First in section
\ref{section-theta} we introduce some notations and recall the
construction from \cite{AKC} of a central extension of the group of
ideles on a smooth projective curve $C$. We show its relation with
the theta-bundle on the Picard variety $\Pic^{g-1}(C)$ of degree
$g-1$ line bundles on $C$. Then in section \ref{secondconstr2} we
give some general construction of a quotient biextension associated
to a bilinear pairing between abelian groups. In section
\ref{tame-Poinc} this construction is applied to the pairing of the
ideles group given by the product of the norms of the tame symbols
over all points of $C$. We show that this defines the Poincar\'e
biextension of the Jacobian of $C$, using results from the previous
sections.

The author is grateful to A. N. Parshin for his help and attention
to this paper, and also thanks L. Breen, D. Osipov, and M. Mazo for
many useful remarks and suggestions.

\section{A central extension of ideles and the theta-bundle}\label{section-theta}

Consider a smooth projective curve $C$ of genus $g$ over a field
$k$. Suppose $K=k(C)$ is the field of rational functions on $C$,
$\A_C=\prod\limits_{x\in C}{'K_x}$ is the ring of ad\`eles on $C$,
and $\A^*_C=\prod\limits_{x\in C}{'K^*_x}$ is the group of ideles.
We put $\Oc=\prod\limits_{x\in C}\hat{\Oc}_x\subset\A_C$ and
$\Oc^*=\prod\limits_{x\in C}\hat{\Oc}^*_x\subset\A_C^*$. There is a
natural surjective homomorphism $\A^*_C\to\Pic(C)$, given by the
formula $\alpha\mapsto [\div(\alpha)]$, where
$\div(\alpha)=\sum\limits_{x\in C}\ord_x(\alpha_x)\cdot x$ and
$\ord_x:K^*_x\to \Z$ is the discrete valuation associated with a
point $x\in C$. The kernel of this homomorphism is equal to the
subgroup $K^*\cdot\Oc^*\subset \A_C^*$. We put
$\deg(\alpha)=\sum\limits_{x\in C}\ord_x(\alpha_x)$. For any two
elements $\alpha_x,\beta_x\in K^*_x$, the tame symbol
$(\alpha_x,\beta_x)_x\in k(x)^*$ is defined by the formula
$$
(\alpha_x,\beta_x)_x=(-1)^{\ord_x(\alpha_x)\ord_x(\beta_x)}
(\alpha_x^{\ord_x(\beta_x)} \beta_x^{-\ord_x(\alpha_x)})(x).
$$

Let us recall some constructions from \cite{AKC}. For any two ideles
$\alpha,\beta\in\A_C^*$, the subspaces $\alpha\Oc\subset\A_C$ and
$\beta\Oc\subset\A_C$ are commensurable, i.e., there exists a
$k$-subspace $L\subset\A_C$ such that $L\subset\alpha\Oc$, $L\subset
\beta\Oc$, and the quotients $(\alpha\Oc)/L$ and $(\beta\Oc)/L$ are
finite dimensional. We put
$(\alpha\Oc|\beta\Oc)=\det_k((\alpha\Oc)/L)^{-1}\otimes\det_k((\beta\Oc)/L)$.
It is easily seen that this does not depend on the choice of $L$ and
$(\alpha\Oc|\beta\Oc)$ is a well defined one-dimensional space over
$k$. The set $\widehat{\A}^*_C=\{(\alpha,r)|\alpha\in\A^*_C,
r\in(\Oc|\alpha\Oc),r\ne 0\}$ has the structure of a group: the
multiplication in $\widehat{\A}^*_C$ is defined by the canonical
isomorphisms $(\Oc|\beta\Oc)\stackrel{\alpha\cdot}\to
(\alpha\Oc|\alpha\beta\Oc)$ and
$(\Oc|\alpha\Oc)\otimes(\alpha\Oc|\alpha\beta\Oc)\to
(\Oc|\alpha\beta\Oc)$ for any ideles $\alpha,\beta\in\A^*_C$. Thus
we get a central extension
$$
1\to k^*\to \widehat{\A}^*_C\to \A^*_C\to 1.\eqno{(*)}
$$
Recall that for any central extension of an abelian group $A$ by an
abelian group $N$
$$
1\to N\to G\stackrel{\pi}\to A\to 1,
$$
the commutator $[g_a,g_b]=g_a g_b g_a^{-1}g_b^{-1}\in N$ depends
only in $a$ and $b$ for any two elements $g_a\in \pi^{-1}(a)$,
$g_b\in \pi^{-1}(b)$; we put $\langle a,b\rangle=[g_a,g_b]$. The
pairing $\langle\cdot,\cdot\rangle$ is skew-symmetric and bilinear.
The first assertion is trivial, while for the second one follows
from the identity
$$
[fg,h]=[f,h][g,h]\Ad(g^{-1}h)(f^{-1})\Ad(hg^{-1})(f)
$$
for all elements $f,g,h\in G$, where $\Ad(g)$ is the conjugation by
$g$. Since the commutator of any two elements in $G$ is central, we
have $\Ad(g^{-1}h)=\Ad(hg^{-1})$ and $[fg,h]=[f,h][g,h]$.

\begin{rmk}
Suppose that the central extension corresponds to the cocycle
$\alpha\in H^2(A,N)$ and let $\bar{\alpha}\colon A\times A\to N$ be
any representative of this cocycle; then $\langle
a,b\rangle=\bar{\alpha}(a,b)\bar{\alpha}(b,a)^{-1}$.
\end{rmk}

The next result was essentially proved in \cite{AKC}.

\begin{theor}\label{AKC-theorem}
For all $\alpha,\beta\in\A^*_C$, the commutator of their liftings to
the group $\widehat{\A}^*_C$ is equal up to sign to the product of
the norms of the tame symbols:
$$
\langle\alpha,\beta\rangle=(-1)^{\deg(\alpha)\deg(\beta)}
\prod_{x\in C}\Nm_{k(x)/k}[(\alpha_x,\beta_x)_x].
$$
\end{theor}

There is a cohomological interpretation of the one-dimensional space
$(\Oc|\alpha\Oc)$, $\alpha\in\A^*_C$. Namely, for any $k$-subspace
$L\subset\A_C$, consider the {\it adelic complex}
$$
\A(L)^{\bullet}\colon 0\to K\oplus L\to \A_C\to 0,
$$
where the differential is given by the formula
$(f,\{f_x\})\mapsto\{f-f_x\}$ for $f\in K$, $\{f_x\}\in L$. Let
$\Df(L)=\det_k H^0(\A(L)^{\bullet})\otimes \det_k
H^1(\A(L)^{\bullet})^{-1}$ be the determinant of cohomology of this
complex. We claim that there is a canonical isomorphism
$$
\Df(\Oc)\otimes(\Oc|\alpha\Oc)\cong \Df(\alpha \Oc).
$$
Indeed, let $L\subset\A_C$ be a $k$-subspace such that $L\subset
\Oc$, $L\subset \alpha\Oc$, and the $k$-spaces $\Oc/L$ and
$(\alpha\Oc)/L$ are finite-dimensional; then the natural embeddings
of complexes $\A(L)^{\bullet}\hookrightarrow\A(\Oc)^{\bullet}$ and
$\A(L)^{\bullet}\hookrightarrow\A(\alpha \Oc)^{\bullet}$ imply the
needed result. In other words,
$(\Oc|\alpha\Oc)=\Hom_k(\Df(\Oc),\Df(\alpha \Oc))$.

This interpretation allows to construct a canonical element
$\widehat{f}\in (\Oc|f\Oc)\backslash\{0\}\subset \widehat{\A}^*_C$
for any $f\in K^*$, using the isomorphism of complexes
$\A(\Oc)^{\bullet}\stackrel{\cdot
f}\longrightarrow\A(f\Oc)^{\bullet}$, which leads to the isomorphism
of one-dimensional $k$-spaces $\Df(\Oc)\to \Df(f\Oc)$. It is easy to
check that the assignment $f\mapsto \widehat{f}$ gives a splitting
of the central extension $(*)$ over the subgroup $K^*\subset
\A^*_C$, i.e, we have $\widehat{f}\cdot\widehat{g}=\widehat{f\cdot
g}$ for all $f,g\in K^*$.

\begin{rmk}
As shown in \cite{AKC}, combining the splitting of $(*)$ over $K^*$
with the formula from Theorem \ref{AKC-theorem}, we get the Weil
reciprocity law on $C$.
\end{rmk}

Further, for any element $u\in \Oc^*\subset \A^*_C$, we have
$\Oc=u\Oc$, hence the space $(\Oc|u\Oc)$ is canonically isomorphic
to $k$. This defines the splitting $\Oc^*\to \widehat{\A}^*_C$,
$u\mapsto\widetilde{u}$ of the extension $(*)$.  We denote the
splittings over $K^*$ and $\Oc^*$ in the different ways because they
do not coincide on the intersection $k^*=K^*\cap \Oc^*\subset
\A^*_C$. Indeed, for a constant $c\in k^*$, we have
$\widehat{c}=c^{\chi(\A(\Oc)^{\bullet})}\cdot\widetilde{c}$. To
compute the Euler characteristic $\chi(\A(\Oc)^{\bullet})$, we give
the following geometrical interpretation of the complexes
$\A(\alpha\Oc)^{\bullet}$, $\alpha\in\A^*_C$.

For any invertible sheaf $\Lc$ on $C$, consider the {\it adelic
complex}
$$
\A(C,\Lc)^{\bullet}\colon 0\to\Lc_{\eta}\oplus\prod_{x\in
C}\hat{\Lc}_x\to \prod_{x\in C}{'\hat{\Lc}_x\otimes_{\hat{\Oc}_x}}
K_x\to 0,
$$
where $\eta$ is the generic point of $C$,
$\hat{\Lc}_x=\Lc_x\otimes_{\Oc_x}\hat{\Oc}_x$ and $\prod'$ is the
adelic product (for more details see \cite{FP} or \cite{Par}). It is
known that there are canonical isomorphisms $H^i(C,\Lc)\cong
H^i(\A(C,\Lc)^{\bullet})$ for $i=0,1$. On the other hand, there is
an equality of complexes
$$
\A(C,\Oc_C(D))^{\bullet}=\A(\alpha\Oc)^{\bullet}
$$
for any idele $\alpha\in\A^*_C$, where $D=-\div(\alpha)$ and for any
open subset $U\subset C$, the group $\Oc_C(D)(U)$ consists of all
functions $f\in K^*$ such that $(\div(f)+D)|_U\ge 0$. Thus there are
canonical isomorphisms $H^i(\A(\alpha\Oc)^{\bullet})\cong
H^i(C,\Oc_C(D))$, $\Df(\alpha\Oc)\cong \det R\Gamma(C,\Oc_C(D))$,
and $(\Oc|\alpha\Oc)\cong\det R\Gamma(C,\Oc_C(D))\otimes\det
R\Gamma(C,\Oc_C)^{-1}$, where $i=0,1$, $\alpha\in\A^*_C$, and
$D=-\div(\alpha)$. In particular, we see that
$\chi(\A(\Oc)^{\bullet})=1-g$, where $g$ is the genus of the curve
$C$ and therefore the splittings of $(*)$ over $K^*$ and $\Oc^*$ do
not coincide in general on the intersection $k^*=K^*\cap \Oc^*$.

\begin{lemma}\label{lemma-cohom-multipl}
\hspace{0cm}

\begin{itemize}
\item[(i)]
For any elements $f\in K^*$ and $\alpha\in\A_C^*$, the following
diagram commutes
$$
\begin{array}{ccc}
(\Oc|\alpha\Oc)&\to&\det R\Gamma(C,\Oc_C(D))\otimes \det R\Gamma(C,\Oc_C)^{-1}\\
\downarrow\lefteqn{\widehat{f}\cdot}&&\downarrow\\
(\Oc|f\alpha\Oc)&\to&\det R\Gamma(C,\Oc_C(D-\div(f)))
\otimes \det R\Gamma(C,\Oc_C)^{-1}\\
\end{array}
$$
where $D=-\div(\alpha)$, the horizontal arrows are canonical
isomorphisms, the first vertical arrow is multiplication on the left
by $\widehat{f}$ in the group $\widehat{\A}^*_C$, and the second
vertical arrow is defined by the canonical isomorphism of invertible
sheaves $\Oc_C(D)\cong\Oc_C(D-(f))$ which is multiplication by $f$
in $K^*$.
\item[(ii)]
For any elements $u\in \Oc^*$ and $\alpha\in\A_C^*$, the following
diagram commutes
$$
\begin{array}{ccc}
(\Oc|\alpha\Oc)&\to&\det R\Gamma(C,\Oc_C(D))\otimes \det R\Gamma(C,\Oc_C)^{-1}\\
\downarrow\lefteqn{\cdot\widetilde{u}}&&\downarrow\lefteqn{id}\\
(\Oc|\alpha u\Oc)&\to&\det R\Gamma(C,\Oc_C(D)))
\otimes \det R\Gamma(C,\Oc_C)^{-1}\\
\end{array}
$$
where $D=-\div(\alpha)$, the horizontal arrows are canonical
isomorphisms, the first vertical arrow is multiplication on the
right by $\widetilde{u}$ in the group $\widehat{\A}^*_C$, and the
second vertical arrow is the identity.
\end{itemize}
\end{lemma}
\begin{proof}
Consider an arbitrary element $r\in(\Oc|\alpha\Oc)\backslash\{0\}$
and the commutative diagram:
$$
\begin{array}{ccc}
\Df(\Oc)&\stackrel{r}\longrightarrow&\Df(\alpha\Oc)\\
\downarrow\lefteqn{f\cdot}& & \downarrow\lefteqn{f\cdot}\\
\Df(f\Oc)&\stackrel{f(r)}\longrightarrow&\Df(f\alpha\Oc).\\
\end{array}
$$
By definition, the composition of the lower triangle in this
diagram, i.e., the diagonal, is equal to $\widehat{f}\cdot r\in
(\Oc|f\alpha\Oc)$. On the other hand, the composition of the upper
triangle corresponds to the identification of $(\Oc|\alpha\Oc)$ with
$(\Oc|f\alpha\Oc)$ via the isomorphism $\det
R\Gamma(C,\Oc_C(D))\stackrel{f\cdot}\longrightarrow \det
R\Gamma(C,\Oc_C(D-\div(f)))$ and this proves $(i)$. The proof of
$(ii)$ is analogous.
\end{proof}

For any integer $n\in\Z$, let $(\A^*_C)^n$ be the set of ideles
$\alpha$ such that $\deg(-\div(\alpha))=n$ and let
$(\widehat{\A}^*_C)^{n}$ be the preimage of $(\A^*_C)^n$ in
$\widehat{\A}^*_C$. Let $\Theta$ be the line bundle on
$\Pic^{g-1}(C)$ whose fiber over an isomorphism class $L$ of degree
$g-1$ line bundles on $C$ is given by $\Theta|_{L}=\det
R\Gamma(C,\Lc)\otimes\det R\Gamma(C,\Oc)^{-1}$, where $\Lc$ is any
representative in $L$. Note that since $\chi(C,\Lc)=0$, this
one-dimensional $k$-space is well defined.

\begin{rmk}
It is known that the line bundle $\Theta$ is isomorphic to the line
bundle associated with the theta-divisor on $\Pic^{g-1}(C)$ (see
\cite{MB}).
\end{rmk}

The next result is a direct consequence of Lemma
\ref{lemma-cohom-multipl}.

\begin{prop}\label{prop-theta-bundle}
There is a well defined action of the group $K^*\cdot \Oc^*$ on the
set $(\widehat{\A}^*_C)^{g-1}$ given by the formula
$(fu)(h)=\widehat{f}\cdot h\cdot\widetilde{u}$ for all $f\in K^*$,
$u\in \Oc^*$, and $h\in (\widehat{\A}^*_C)^{g-1}$; this action
commutes with the natural action of $K^*\cdot \Oc^*$ on
$(\A^*_C)^{g-1}$. Moreover, there is a canonical isomorphism of
$k^*$-torsors on $\Pic^{g-1}(C)$
$$
K^*\backslash(\widehat{\A}^*_C)^{g-1}/\Oc^*\cong
\Theta\backslash\{0\},
$$
where we identify $K^*\backslash(\A^*_C)^{g-1}/\Oc^*$ with
$\Pic^{g-1}(C)$ via the map $\alpha\mapsto -\div(\alpha)$,
$\alpha\in\A^*_C$.
\end{prop}

\section{Construction of a quotient biextension}\label{secondconstr2}

For all groups below, we write the group law in the multiplicative
way. See more details on biextensions in \cite{SGA} and \cite{Bre}.
Let $A,A',N$ be abelian groups, $B,B'\subset A$, $C,C'\subset A$ be
subgroups, and let $\langle\cdot,\cdot\rangle:A\times A\to N$ be a
bilinear pairing such that $\langle B,B'\rangle=1$, $\langle
C,C'\rangle=1$, $\langle B\cap C,A'\rangle=1$, and $\langle A,B'\cap
C'\rangle=1$. Let $T$ be the trivial biextension of $(A,A')$ by $N$.
By $T|_{(a,a')}$ denote the fiber of $T$ over $(a,a')\in A\times
A'$. For all elements $a\in A$, $bc\in B\cdot C$, $a'\in A'$, and
$b'c'\in B'\cdot C'$, consider the isomorphism $T|_{(a,a')}\to
T|_{(abc,a'b'c')}$ that is equal to multiplication by $\langle
a',c\rangle\langle b',a\rangle\langle b',c\rangle\in N$. It is
readily seen that the last expression does not depend on the
decompositions $bc$ and $b'c'$ and it can be checked that this
defines an action of the group $(B\cdot C)\times(B'\cdot C')$ on
$T$, which commutes with the natural action on $A\times A'$.
Moreover, this action commutes with the biextension structure on $T$
and we get the quotient biextension $P=T/((B\cdot C)\times(B'\cdot
C'))$ of $(A/(B\cdot C),A'/(B'\cdot C'))$ by $N$.

\begin{rmk}\label{rmk-coin-actions}
Given a central extension $1\to N\to G\to A\to 1$, one has a
canonically trivial biextension $\Lambda(G)=m^*G\wedge
p_1^*G^{-1}\wedge p_2^*G^{-1}$ of $(A,A)$ by $N$, where $p_1$,
$p_2$, $m$, and $\wedge$ denote, respectively, projection on the
first multiple, projection on the second multiple, multiplication in
the group $A$, and product in the category of $N$-torsors on
$A\times A$. If the extension splits over subgroups $B,C\subset A$,
then the commutator pairing $\langle\cdot,\cdot\rangle$ satisfies
the above condition with $A'=A$, $B'=B$, and $C=C'$ if we suppose
also that $\langle B\cap C,A\rangle=1$. Let the assignments
$b\mapsto \widehat{b}$ and $c\mapsto \widetilde{c}$ be splittings of
the given central extension over $B$ and $C$, respectively; then
there is an action of $B\times C$ on $G$ given by the formula
$g\mapsto\hat{b}g\widetilde{c}$. This naturally induces the action
of $(B\times C)\times (B\times C)$ on $\Lambda(G)$. An explicit
calculation shows that this action factors through $(B\cdot C)\times
(B\cdot C)$ and coincides with the one defined above.
\end{rmk}

Let $D$ and $D'$ be two abelian groups. Recall that for any
biextension $P$ of $(D,D')$ by $N$ and for any integer $m\in\Z$,
$m\ge 1$, one defines the {\it Weil pairing} $\phi_m:A_m\times
A_m\to N_m$ in the following way: for any $(d,d')\in D_m\times
D'_m$, the element $\phi_m(d,d')$ is equal to the composition given
by the diagram of $N$-torsors
$$
\begin{array}{ccc}
P^{\wedge m}|_{(d,d')}&\longrightarrow&P|_{(d^m,d')}\\
\uparrow&&\downarrow\\
P|_{(d,(d')^m)}&\longleftarrow& P|_{(1,1)},\\
\end{array}
$$
where the arrows are natural isomorphisms of $N$-torsor defined by
the biextension structure on $P$.

\begin{prop}\label{Weil-pair-commut}
Let $A,A',B,B',C,C',\langle\cdot,\cdot\rangle,T,P$ be as in the
beginning of this section and let $a\in A,a'\in A',b\in B,b'\in
B,c\in C,c'\in C'$ be such that $a^m=bc$ and $(a')^m=b'c'$; then we
have
$$
\phi_m(\bar{a},\bar{a}')=\langle b',a\rangle\langle a',
c\rangle^{-1},
$$
where $\bar{a}\in A/(B\cdot C)$ and $\bar{a}'\in A'/(B'\cdot C')$
are the classes corresponding to $a$ and $a'$, respectively.
\end{prop}
\begin{proof}
By construction, the pull-back of the biextension $P$ from
$A/(B\cdot C)\times A'/(B'\cdot C')$ to $A\times A'$ is isomorphic
to the trivial biextension $T$. Therefore the pull-back of the
diagram defining the Weil pairing $\phi_m(\bar{a},\bar{a}')$ is the
diagram
$$
\begin{array}{ccc}
T^{\wedge
m}|_{(a,a')}&\stackrel{id}\longrightarrow&T|_{(bc,a')}\\
\uparrow\lefteqn{id}&&\downarrow\lefteqn{\langle
a',c \rangle^{-1}}\\
T|_{(a,b'c')}&\stackrel{\langle b',a\rangle}\longleftarrow& T|_{(1,1)}.\\
\end{array}
$$
This concludes the proof.
\end{proof}

\section{Tame symbols and the Poincar\'e
biextension}\label{tame-Poinc}

As before, let $C$ be a smooth projective curve of genus $g$ over a
field $k$. Let us recall a construction of the Poincar\'e
biextension $\Pc$ of $(\Pic^0(C),\Pic^0(C))$ by $k^*$ (see
\cite{Del} and \cite{MB}). For all isomorphism classes $L,M$ of
degree zero line bundles on $C$, we put $\Pc|_{(L,M)}=(\Lc,\M)$,
where
$$
(\Lc,\M)=(\det R\Gamma(C,\Lc\otimes\M)\otimes \det
R\Gamma(C,\Lc)^{-1}\otimes \det R\Gamma(C,\M)^{-1}\otimes\det
R\Gamma(C,\Oc_C))\backslash\{0\}
$$
and $\Lc$ and $\M$ are any representatives from $L$ and $M$,
respectively. Since
$\chi(C,\Lc\otimes\M)=\chi(C,\Lc)=\chi(C,\M)=1-g$, this
one-dimensional $k$-space is well defined and $\Pc$ is a
$k^*$-torsor on $\Pic^0(C)\times\Pic^0(C)$. To define a biextension
structure on $\Pc$, consider the pull-back $p^*\Pc$ of $\Pc$ with
respect to the natural map $p:\Pic^0(C)\times\Div^0(C)\to
\Pic^0(C)\times\Pic^0(C)$ given by the formula $(L,D)\mapsto
(L,[\Oc_C(D)])$. There is a canonical isomorphism
$\varphi:p^*\Pc\cong \Pc'$ of $k^*$-torsors on
$\Pic^0(C)\times\Div^0(C)$, where $\Pc'$ is the biextension of
$(\Pic^0(C),\Div^0(C))$ by $k^*$ defined by the formula
$\Pc'|_{(L,D)}=(\bigotimes\limits_{x\in
C}\Lc|_x^{\otimes\ord_x(D)})\backslash\{0\}$, where $\Lc$ is any
representative from the class $L$. Thus $\varphi$ induces a
biextension structure on $p^*\Pc$ and it turns out that it descends
to a biextension structure on $\Pc$.

Now we put $A=A'=(\A^*_C)^0$, $B=B'=K^*$, $C=C'=\Oc^*$, $N=k^*$, and
$\langle\alpha,\beta\rangle=\prod\limits_{x\in C}{\rm
Nm}_{k(x)/k}[(\alpha_x,\beta_x)_x]$. It is readily seen that the
conditions from the beginning of section \ref{secondconstr2} are
satisfied, hence we get a biextension $P$ of $(\Pic^0(C),\Pic^0(C))$
by $k^*$.

\begin{theor}\label{thm-quot-Poinc}
The biextension $P$ is canonically isomorphic to the Poincar\'e
biextension $\Pc$ of $(\Pic^0(C),\Pic^0(C))$ by $k^*$.
\end{theor}
\begin{proof}
Let $\pi':(\A^*_C)^0\times(\A^*_C)^0\to \Pic^0(C)\times\Div^0(C)$ be
the homomorphism given by the formula $(\alpha,\beta)\mapsto
([\Oc_C(-\div(\alpha))],-\div(\beta))$ and let
$\pi=p\circ\pi':(\A^*_C)^0\times(\A^*_C)^0\to
\Pic^0(C)\times\Pic^0(C)$.

It follows from section \ref{section-theta} that there is a
canonical isomorphism
$\psi:\Lambda((\widehat{\A}^*_C)^0)\to\pi^*\Pc$ of $k^*$-torsors
over $(\A^*_C)^0\times(\A^*_C)^0$. Combining Theorem
\ref{AKC-theorem}, Lemma \ref{lemma-cohom-multipl}, and Remark
\ref{rmk-coin-actions}, we see that the natural action of
$(K^*\cdot\Oc^*)\times(K^*\cdot\Oc^*)$ on $\pi^*\Pc$ commutes via
$\psi$ with the action of $(K^*\cdot\Oc^*)\times(K^*\cdot\Oc^*)$ on
the canonically trivial biextension $\Lambda((\widehat{\A}^*_C)^0)$
described in the beginning of section \ref{secondconstr2} and
defining the biextension $P$. Therefore $P$ is isomorphic to $\Pc$
as a $k^*$-torsor on $\Pic^0(C)\times\Pic^0(C)$ and it remains to
check that the canonical isomorphism $\psi$ commutes with the
biextension structures.

Note that the pull-back $(\pi')^*\Pc'$ has a canonical
trivialization given by the assignment
$$
(\alpha,\beta)\mapsto \mbox{$\bigotimes\limits_{x\in C}
\bar{\alpha}_x^{\otimes(-\ord_x(\beta_x))}\in\bigotimes\limits_{x\in
C} (\alpha_x\Oc_x/\m_x\alpha_x\Oc_x)^{\otimes(-\ord_x(\beta_x))}=
(\pi')^*\Pc'|_{(\alpha,\beta)}$}.
$$
Thus it suffices to check that the composition $(\pi')^*\varphi\circ
\psi:\Lambda((\widehat{\A}^*_C)^0)\to (\pi')^*\Pc'$ sends one
trivialization to the other.

Let us recall the explicit form of the isomorphism $\varphi$. Take a
pair $([\Lc],D)\in\Pic^0(C)\times\Div^0(C)$. First, suppose that
$D\ge 0$; then the exact sequences of sheaves
$$
0\to\Lc\to\Lc(D)\to\Lc(D)|_D\to 0,
$$
$$
0\to\Oc_'\to\Oc_'(D)\to\Oc_'(D)|_D\to 0
$$
lead to the isomorphism
$\mu:(\Lc,\Oc_C(D))\to\det_k(\Lc(D)|_D)\otimes
\det_k(\Oc_'(D)|_D)^{-1}$. Further, by induction on the degree of
$D$, one establishes a canonical isomorphism
$$
\mbox{$\nu\colon\Hom_k(\det_k(\Oc_C(D)|_D),\det_k(\Lc(D)|_D))\cong
\bigotimes\limits_{x\in C}\Lc|_x^{\otimes\ord_x(D)}$}.
$$
The isomorphism $\varphi$ equals the composition $\nu\circ\mu$.

Suppose that $\{s_x\}\in\prod\limits_{x\in C}\Lc_x$ is a collection
of local sections such that $\bar{s}_x\ne 0$ for all $x\in C$, where
$\bar{s}_x\in\Lc|_x$ is the value at a point $x$ of a section
$s_x\in\Lc_x$. Then the determinant of the isomorphism
$\bigotimes\limits_{x\in |D|}s_x\colon\Oc_C(D)|_D\to\Lc(D)|_D$ is
mapped under $\nu$ to the product $\bigotimes\limits_{x\in
C}\bar{s}_x^{\otimes\ord_x(D)}$, where $|D|$ is the support of the
divisor $D$.

Now consider the pull-backs of $\mu$ and $\nu$ with respect to
$\pi'$. Let $(\alpha,\beta)\in(\A^*_C)^0\times (\A^*_C)^0$ be such
that $\pi'(\alpha,\beta)=([\Lc],D)$. We may assume that
$\Lc=\Oc_C(-\div(\alpha))$. Then the map $(\pi')^*\mu\circ \psi$ is
the natural isomorphism
$$
\Lambda((\widehat{\A}^*_C)^0)|_{(\alpha,\beta)}=(\Oc|\alpha\beta\Oc)\otimes
(\Oc|\alpha\Oc)^{-1}\otimes(\Oc|\beta\Oc)^{-1}\to
(\alpha\Oc|\alpha\beta\Oc)\otimes(\Oc|\beta\Oc)^{-1}
$$
that follows from the exact sequences of complexes
$$
0\to\A(\alpha\Oc)\to\A(\alpha\beta\Oc)\to
(\alpha\Oc|\alpha\beta\Oc)\to 0,
$$
$$
0\to\A(\Oc)\to\A(\beta\Oc)\to (\Oc|\beta\Oc)\to 0.
$$
Therefore the isomorphism $(\pi')^*\mu\circ\psi$ takes the canonical
element in $\Lambda((\widehat{\A}^*_C)^0)|_{(\alpha,\beta)}$ to the
element
$\det(\alpha)\in\Hom_k((\Oc|\beta\Oc),(\alpha\Oc|\alpha\beta\Oc))$
that equals to the determinant of the isomorphism
$(\Oc|\beta\Oc)\stackrel{\alpha\cdot}\to
(\alpha\Oc|\alpha\beta\Oc)$. Further, the idele $\alpha$ defines a
collection of local sections $\{\alpha_x\}\in\prod\limits_{x\in
C}\Lc_x$, hence $(\pi')^*\nu(\det(\alpha))=\bigotimes\limits_{x\in
C} \bar{\alpha}_x^{\otimes(-\ord_x(\beta_x))}$. Thus we have treated
the case when the divisor $D$ is effective.

One considers the case when $E=-D\ge 0$ in the same way, using the
exact sequences of sheaves
$$
0\to\Lc(-E)\to\Lc\to\Lc|_E\to 0,
$$
$$
0\to\Oc_C(-E)\to\Oc_C\to\Oc_C|_E\to 0.
$$
The case of an arbitrary divisor $D-E$, where $D,E\ge 0$, can be
reduced to these two cases, using the embeddings of sheaves
$\Lc(-E)\subset\Lc$ and $\Lc(-E)\subset\Lc(D-E)$ (respectively,
$\Oc_C(-E)\subset\O$ and $\Oc_C(-E)\subset\Oc_C(D-E)$), whose
pull-back with respect to $\pi'$ will correspond to the choice of a
common $k$-subspace in the commensurable spaces $\alpha\Oc$ and
$\alpha\beta\Oc$ (respectively, $\Oc$ and $\beta\Oc$) when defining
the one-dimensional space $(\alpha\Oc|\alpha\beta\Oc)$
(respectively, $(\Oc|\beta\Oc)$).
\end{proof}

\begin{rmk}
One can also descend a symmetric structure from the trivial
biextension $T$ of $((\A^*_C)^0,(\A^*_C)^0)$ to the biextension $P$
and check this coincides with the natural symmetric structure on
$\Pc$.
\end{rmk}

Recall that for a natural number $m$ prime to ${\rm char}(k)$, the
Weil pairing $\phi_m:\Pic^0(C)_m\times\Pic^0(C)_m\to \mu_m$ is the
Weil pairing in the above sense associated with the Poincar\'e
biextension $\Pc$ (see \cite{MB}). Combining Proposition
\ref{Weil-pair-commut} with Theorem \ref{thm-quot-Poinc}, we get the
following adelic formula for the Weil pairing.

\begin{corol}\label{old-work}
Let $\alpha,\alpha'\in\A^*_C$ be two ideles such that $\alpha^m=fu$
and $(\alpha')^m=f'u'$, where $f,f'\in K^*$, $u,u'\in \Oc^*$, and
let $\Lc=\Oc_C(-\div(\alpha))$, $\M=\Oc_C(-\div(\alpha'))$; then we
have
$$
\phi_m(\Lc,\M)=\prod_{x\in
C}\Nm_{k(x)/k}[(f,\alpha'_x)_x(\alpha_x,u'_x)^{-1}_x].
$$
\end{corol}

\begin{rmk}
If the divisors $D=-\div(\alpha)$ and $D'=-\div(\alpha')$ do not
intersect, then we have $\phi_m(\Lc,\M)=f'(D)\cdot f^{-1}(D')$. The
equivalence of this definition of the Weil pairing with usual one
was first shown by Howe in \cite{How}.
\end{rmk}

\begin{rmk}
Suppose that the ground field $k$ is algebraically closed; then the
group $\hat{\Oc}_x$ is $m$-divisible for any closed point $x\in C$
and any integer $m$ prime to ${\rm char}(k)$. Therefore, given the
divisors $D$ and $D'$, one may choose the ideles $\alpha$ and
$\alpha'$ such that $\alpha^m=f$ and $(\alpha')^m=f'$, where
$f,f'\in K^*$. Then $\phi_m(\Lc,\M)=\prod\limits_{x\in
C}\Nm_{k(x)/k}[(\alpha_x,\alpha'_x)^m_x]$. The coincidence of this
formula with the definition of the Weil pairing via biextensions was
also directly explained by Mazo in \cite{Maz}.

\end{rmk}

\end{document}